\documentclass[11pt,reqno]{amsart}
\usepackage{amssymb,latexsym}
\usepackage{graphicx}
\usepackage{fancyhdr}

\theoremstyle{plain}
\numberwithin{equation}{section}
\newtheorem{thm}{Theorem}[section]

\newtheorem{coll}[thm]{Corollary}

\begin{document}

\setcounter{page}{1}

\title[An explicit formula generating the non-Fibonacci numbers]{An explicit formula generating the non-Fibonacci numbers}
\author{Bakir Farhi}
\address{Department of Mathematics \\
University of B\'ejaia \\
B\'ejaia \\
Algeria
}
\email{bakir.farhi@gmail.com}

\begin{abstract}
We show among others that the formula:
$$
\left\lfloor n + \log_{\Phi}\left\{\sqrt{5}\left(\log_{\Phi}(\sqrt{5}n) + n\right) -5 + \frac{3}{n}\right\} - 2 \right\rfloor ~~~~ (n \geq 2) ,
$$
(where $\Phi$ denotes the golden ratio and $\lfloor \rfloor$ denotes the integer part) generates the non-Fibonacci numbers.
\end{abstract}

\maketitle

\section{Introduction and main result}
Two sequences of natural numbers are said to be complementary if they are disjoint and their union is the entire set $\mathbb{N}$ of nonnegative integers. Given a sequence of nonnegative integers, it is an important problem to find an explicit formula for the sequence complementing them in $\mathbb{N}$. For example, it is well known that the sequence complementing the perfect square positive integers is generated by the formula $\lfloor n + \sqrt{n} + \frac{1}{2}\rfloor$ $(n \in \mathbb{N})$ and the sequence complementing the triangular numbers (i.e., the integers having the form $\frac{n (n + 1)}{2}$, $n \in \mathbb{N}$) is generated by the formula $\lfloor n + \sqrt{2 n} + \frac{1}{2}\rfloor$ $(n \geq 1)$. On this topic we can consult the article \cite{lm} of Lambekand and Moser. For a brief introduction, we can consult the book \cite{hons2} of Honsberger or the chapter 1 of the book \cite{honsberger} of the same author. 

In this paper, we establish a general theorem which gives the formula generating the complement (in $\mathbb{N}$) of a given sequence of integers. Then we apply it to obtain the complementary sequences of some types of sequences including the Fibonacci sequence. So we obtain an explicit formula for the $n$\textsuperscript{th} non-fibonacci number. Actually, Gould \cite{gould} already obtained an approximate formula for the $n$\textsuperscript{th} non-Fibonacci number (noted $g_n$). His formula is:
$$
g_n = n + F(n + F(n + F(n))) ,
$$
where $F(n) = \lfloor \log_{\Phi}n + \frac{1}{2}\log_{\Phi}5 - 1\rfloor$.
But the inconvenient of Gould's formula is that it is quite complicated and the purpose of this paper is to obtain an easy formula depending only on $n$. In addition, our approach is somewhat different from that of Gould. Our main result is the following:
\begin{thm}\label{t1}
Let ${(u_n)}_{n \in \mathbb{N}}$ be an increasing sequence of integers and $\varphi : [0 , + \infty[ \rightarrow \mathbb{R}$ be a continued function which increases and tends to $+ \infty$ when $x$ tends to $+ \infty$. Suppose that $\varphi$ satisfies for all $n \in \mathbb{N}$:
\begin{equation}\label{eq1}
u_n - n < \varphi(n) \leq u_n - n + 1 \tag{$I$}
\end{equation}
Then the formula ${(\lfloor n + \varphi^{-1}(n)\rfloor)}_{n \geq u_0 + 1}$ generates the complement of $\{u_n , n \in \mathbb{N}\}$ in $[u_0 , + \infty[ \cap \mathbb{Z}$.
\end{thm}

\noindent {\bf Proof.}\\
First remark that the hypothesis \eqref{eq1} of the theorem gives $\varphi(0) \leq u_0 + 1$. Consequently, since $\varphi$ is continuous, increasing and tends to $+ \infty$ when $x$ tends to $+ \infty$, the set of arrival of $\varphi$ contains the interval $[u_0 + 1 , + \infty[$. It follows that $\varphi^{-1}$ is defined at least in the interval $[u_0 + 1 , + \infty[$. Also, because ${(u_n)}_{n \in \mathbb{N}}$ is an increasing sequence of integers, we have for all $n \in \mathbb{N}$: $u_n \geq u_0 + n$. In particular $\varphi^{-1}$ can be applied to all real number greater than or equal to $u_n - n + 1$ $(n \in \mathbb{N})$.\\
$\bullet$ Now, let $N \in [u_0 , + \infty[ \cap \mathbb{Z}$ which is not a term of ${(u_n)}_n$ and let us show that $N$ is a term of the sequence ${(\lfloor k + \varphi^{-1}(k)\rfloor)}_{k \geq u_0 + 1}$. Since $N \geq u_0$ then $N$ lies between two consecutive terms of ${(u_n)}_n$. Let $n \in \mathbb{N}$ such that:
$$
u_n < N < u_{n + 1} .
$$
since $N$ is an integer, we have also:
$$
u_n + 1 \leq N \leq u_{n + 1} - 1 .
$$
Hence
$$
u_n - n + 1 \leq N - n \leq u_{n + 1} - (n + 1) .
$$
By applying $\varphi^{-1}$ (which is increasing because $\varphi$ is increasing), it follows that:
$$
\varphi^{-1}(u_n - n + 1) \leq \varphi^{-1}(N - n) \leq \varphi^{-1}(u_{n + 1} - (n + 1)) .
$$
But according to the hypothesis \eqref{eq1}, we have $\varphi^{-1}(u_{n + 1} - (n + 1)) < n + 1$ and $\varphi^{-1}(u_n - n + 1) \geq n$. So
$$
n \leq \varphi^{-1}(N - n) < n + 1 ,
$$
which implies that:
$$
\lfloor \varphi^{-1}(N - n)\rfloor = n .
$$
Finally, we conclude that:
$$
N = (N - n) + n = (N - n) + \lfloor \varphi^{-1}(N - n)\rfloor = \lfloor (N - n) + \varphi^{-1}(N - n)\rfloor ,
$$
which implies that $N$ is generated by the formula $\lfloor k + \varphi^{-1}(k)\rfloor$ $(k \geq u_0 + 1)$. \\
$\bullet$ Conversely, let $N$ be a term of the sequence ${(\lfloor k + \varphi^{-1}(k)\rfloor)}_{k \geq u_0 + 1}$ and let us show that $N$ is not a term of ${(u_k)}_{k \in \mathbb{N}}$.\\
Let $n \geq u_0 + 1$ be fixed such that $N = \lfloor n + \varphi^{-1}(n)\rfloor$. So we have:
$$
N \leq n + \varphi^{-1}(n) < N + 1 .
$$
By subtracting $n$ and then applying $\varphi$ (which is increasing), we get:
$$
\varphi(N - n) \leq n < \varphi(N - n + 1) .
$$
But according to the hypothesis \eqref{eq1}, we have: $\varphi(N - n + 1) \leq u_{N - n + 1} - (N - n + 1) + 1 = u_{N - n + 1} - N + n$ and $\varphi(N - n) > u_{N - n} - (N - n) = u_{N -n} - N + n$. Using this, we get:
$$
u_{N - n} - N + n < n < u_{N - n + 1} - N + n ,
$$
which is equivalent to:
$$
u_{N - n} < N < u_{N - n + 1} .
$$
So $N$ lies between two consecutive terms of ${(u_k)}_k$. Hence $N$ cannot be a term of ${(u_k)}_k$.\\
This completes the proof of the theorem.\hfill$\blacksquare$

\section{Applications}
\subsection{The Complement of the sequence ${(n^r)}_{n \in \mathbb{N}}$, $r \in \mathbb{N} , r \geq 2$}
We have the following:
\begin{coll}
Let $r \geq 2$ be an integer. The formula ${\left(\left\lfloor n + \sqrt[r]{n + \sqrt[r]{n}}\right\rfloor\right)}_{n \geq 1}$ generates the positive integers which are not $r$\textsuperscript{th} powers.
\end{coll}

\noindent {\bf Proof.}\\
We apply Theorem \ref{t1} for $u_n = n^r$ $(n \in \mathbb{N})$ and $\psi(x) := \varphi^{-1}(x) = \sqrt[r]{x + \sqrt[r]{x}}$ $(x \in [0 , + \infty[)$ which is a continuous and increasing function and tends to $+ \infty$ when $x$ tends to $+ \infty$. To verify the hypothesis \eqref{eq1} of Theorem \ref{t1}, it is equivalent to verify that
$$
\psi(n^r - n) < n ~~\text{and}~~ \psi(n^r - n + 1) \geq n ~~~~ (\forall n \in \mathbb{N}) .
$$
In odrer to show that $\psi(n^r - n) < n$, it suffices to bound from above $\sqrt[r]{n^r - n}$ by $n$ and in order to show that $\psi(n^r - n + 1) \geq n$ it suffices to bound from bellow $\sqrt[r]{n^r - n + 1}$ by $(n - 1)$. The corollary follows.\hfill$\blacksquare$

\medskip

\noindent {\bf Remark.}
For the sequences of perfect squares and perfect cubes, we have other formulas more sample than the previous one complementing them. Indeed, using Theorem \ref{t1}, we can show that the formula ${(\lfloor n + \sqrt{n} + \frac{1}{2}\rfloor)}_{n \geq 1}$ generates the positive integers which are not perfect squares and the formula ${\left(\left\lfloor n + \sqrt[3]{n} + \frac{1}{3 \sqrt[3]{n + 1}}\right\rfloor\right)}_{n \geq 1}$ generates the positive integers which are not perfect cubes.

\subsection{The complement of the sequence ${(a^n)}_{n \in \mathbb{N}}$, $a \in \mathbb{N}$, $a \geq 2$}
We have the following:
\begin{coll}
Let $a \geq 2$ be an integer. The formula $\left\lfloor n + \log_a\left(n + \log_a(n)\right)\right\rfloor$ $(n \geq 1)$ generates the positive integers which are not powers of $a$.
\end{coll}

\noindent{\bf Proof.}\\
We apply Theorem \ref{t1} for $u_n = a^n$ $(n \in \mathbb{N})$ and $\psi(x) := \varphi^{-1}(x) = \log_a\left(x + \log_a(x)\right)$ $(x \in [1 , + \infty[)$ which is a continuous and increasing function and tends to $+ \infty$ when $x$ tends to $+ \infty$.\\ To verify the hypothesis \eqref{eq1} of Theorem \ref{t1}, we have to verify that:
$$
\psi(a^n - a) < n ~~\text{and}~~ \psi(a^n - n + 1) \geq n ~~~~ (\forall n \in \mathbb{N}) .
$$
Those inequalities easily follow from the trivial upper bound $\log_a(a^n - n) < n$ and the trivial lower bound $\log_a(a^n - n + 1) > n - 1$. The corollary follows.\hfill$\blacksquare$ 

\subsection{The complement of the Fibonacci sequence}
The Fibonacci sequence, noted ${(F_n)}_{n \in \mathbb{N}}$, is defined by:
$$
\left\{
\begin{array}{l}
F_0 = 0 ~,~ F_1 = 1 \\
F_{n + 2} = F_n + F_{n + 1} ~~~~ (\forall n \in \mathbb{N})
\end{array}
\right..
$$
The Fibonacci numbers are simply the terms of ${(F_n)}_n$. The complement of ${(F_n)}_n$ is given by the following:
\begin{coll}
The formula
$$
\left\lfloor n + \log_{\Phi}\left\{\sqrt{5}\left(\log_{\Phi}(\sqrt{5}n) + n\right) -5 + \frac{3}{n}\right\} - 2 \right\rfloor ~~~~ (n \geq 2)
$$
generates the numbers which are not Fibonacci numbers.
\end{coll}

\noindent {\bf Proof.}\\
We apply Theorem \ref{t1} for $u_n = F_{n + 2}$ $(n \in \mathbb{N})$ and
$$
\psi(x) := \varphi^{-1}(x) = \log_{\Phi}\left\{\sqrt{5}\left(\log_{\Phi}(\sqrt{5} x) + x\right) -5 + \frac{3}{x}\right\} - 2 ,
$$
which is a continuous and increasing function on $[2 , + \infty[$ and tends to $+ \infty$ when $x$ tends to $+ \infty$.\\
To verify the hypothesis \eqref{eq1} of Theorem \ref{t1}, we have to verify that:
$$
\psi(F_{n + 2} - n) < n ~~\text{and}~~ \psi(F_{n + 2} - n + 1) \geq n ~~~~ (\forall n \in \mathbb{N}) .
$$
To do so, we verify those inequalities for the small values of $n$ ($n \leq 10$) and we use Binet's formula (see for example \cite{honsberger}, chapter 8):
$$
F_n = \frac{1}{\sqrt{5}} \left(\Phi^n - \overline{\Phi}^n\right)
$$
(where $\overline{\Phi} = \frac{1 - \sqrt{5}}{2} = -\frac{1}{\Phi}$) to verify them for the large values of $n$ ($n > 10$). Let us prove the above inequalities for the large values of $n$. Using Binet's formula, the calculations give:
\begin{multline*}
\sqrt{5}\left(\log_{\Phi}(\sqrt{5} (F_{n + 2} - n)) + F_{n + 2} - n\right) - 5 + \frac{3}{F_{n + 2} - n} ~=~ \Phi^{n + 2} - \overline{\Phi}^{n + 2} \\ + 2 \sqrt{5} - 5 + \frac{3}{F_{n + 2} - n} + \sqrt{5} \log_{\Phi}\left\{1 - (- \overline{\Phi}^2)^{n + 2} - \sqrt{5} n (- \overline{\Phi})^{n + 2}\right\} .
\end{multline*}
Because $2 \sqrt{5} - 5 < 0$ and the quantity
$$
- \overline{\Phi}^{n + 2} + \frac{3}{F_{n + 2} - n} + \sqrt{5} \log_{\Phi}\left\{1 - (- \overline{\Phi}^2)^{n + 2} - \sqrt{5} n (- \overline{\Phi})^{n + 2}\right\}
$$
tends to $0$ as $n$ tends to infinity, we have:
$$
\sqrt{5}\left(\log_{\Phi}(\sqrt{5} (F_{n + 2} - n)) + F_{n + 2} - n\right) - 5 + \frac{3}{F_{n + 2} - n} < \Phi^{n + 2}
$$
for $n$ sufficiently large (in practice $n > 10$ suffices). This gives $\psi(F_{n + 2} - n) < n$, as required.\\
Similarly, using Binet's Formula, the calculations give:
\begin{multline*}
\sqrt{5}\left(\log_{\Phi}(\sqrt{5} (F_{n + 2} - n + 1)) + F_{n + 2} - n + 1\right) - 5 + \frac{3}{F_{n + 2} - n + 1} ~=\\ \Phi^{n + 2} - \overline{\Phi}^{n + 2} + 3 \sqrt{5} - 5 + \frac{3}{F_{n + 2} - n + 1} + \sqrt{5} \log_{\Phi}\left\{1 - (- \overline{\Phi}^2)^{n + 2} - \sqrt{5} (n - 1) (- \overline{\Phi})^{n + 2}\right\} .
\end{multline*}
Because $3 \sqrt{5} - 5 > 0$ and the quantity
$$
- \overline{\Phi}^{n + 2} + \frac{3}{F_{n + 2} - n + 1} + \sqrt{5} \log_{\Phi}\left\{1 - (- \overline{\Phi}^2)^{n + 2} - \sqrt{5} (n - 1) (- \overline{\Phi})^{n + 2}\right\}
$$
tends to $0$ as $n$ tends to infinity then for $n$ sufficiently large ($n > 10$ suffices) we have:
$$
\sqrt{5}\left(\log_{\Phi}(\sqrt{5} (F_{n + 2} - n + 1)) + F_{n + 2} - n + 1\right) - 5 + \frac{3}{F_{n + 2} - n + 1} > \Phi^{n + 2} .
$$
This gives $\psi(F_{n + 2} - n + 1) > n$ (for $n > 10$). The proof is complete.\hfill$\blacksquare$

\newpage

\noindent {\bf Remark on the sequences with several indices}

\noindent We don't know how to generalize Theorem \ref{t1} for the sequences of several indices although there exist some theorems of complement of sequences with several indices. The more famous is perhaps Legendre's theorem (see e.g. \cite{dav}) which states that the sequence with three indices ${(n^2 + m^2 + k^2)}_{n , m , k \in \mathbb{N}}$ has for complement (in $\mathbb{N}$) the sequence with two indices ${(4^h (8 \ell + 7))}_{h , \ell \in \mathbb{N}}$.\\
Note also that if we are able to generalize Theorem \ref{t1} for sequences with two indices then we can obtain a formula generating prime numbers, because it is obvious that the sample formula ${((n + 2)(m + 2))}_{n , m \in \mathbb{N}}$ generates the composite numbers.

\end{document}